\documentclass[10pt,twoside]{article}
\usepackage{graphicx}
\usepackage{amsmath}
\usepackage{Latex-document}
\def\R{{{{\rm l} \kern -.15em {\rm R}}}}
\def\Z{{{{\rm Z} \kern -.35em {\rm Z}}}}
\def\e{\varepsilon}
\def\dis{\displaystyle}
\def\cT{{\cal T}}
\def\wh{\widehat}
\def\n{\noindent}

\title{\bf  High Dimensional Finite Elements for\vskip -2mm Elliptic Problems with Multiple
Scales\vskip -2mm and Stochastic Data\thanks{Research performed in
the network ``Homogenization and Multiple Scales 2000'' (HMS 2000)
of the EC and supported by the Swiss Federal government under
grant No. BBW 01.0025-1 and by the Swiss National Fund under grant
Number SNF 21-58754.99.}\vskip 6mm}

\author{C. Schwab\vspace*{-0.5cm}\thanks{Seminar for Applied Mathematics,
ETH Z\"urich, 8092 Z\"urich, Switzerland. E-mail: schwab@sam.math.ethz.ch}}
\date{\vspace{-8mm}}

\begin{document}
\maketitle

\thispagestyle{first} \setcounter{page}{727}

\begin{abstract}\vskip 3mm
Multiple scale homogenization problems are reduced to single scale problems in higher dimension.
It is shown that sparse tensor product Finite Element Methods (FEM) allow the numerical solution
in complexity independent of the dimension and of the length scale.
Problems with stochastic input data are reformulated as
high dimensional deterministic problems for the statistical moments of the random solution.
Sparse tensor product FEM give a deterministic solution algorithm of log-linear complexity
for statistical moments.

\vskip 4.5mm


\end{abstract}

\vskip 12mm

\section{Introduction} \label{section4 1}\setzero
\vskip-5mm \hspace{5mm}

The numerical solution of elliptic problems with multiple scales in a bounded
domain $D$
can be achieved either by analytic homogenization \cite{BLP78}, \cite{ST01}
through asymptotic analysis or by specially designed Finite Element spaces to capture
the fine scales of the problem \cite{MBS2000}, \cite{MS2002}.
In asymptotic theory of homogenization, the fine scale of the solution is averaged and
lost in the homogenized limit. Fine scale information can be recovered from so called
correctors which must be calculated separately.
An alternative which we pursue here is to ``unfold'' the homogenization problem into a
single scale problem in high dimension - so that homogenized and fine scale behaviour
are still coupled \cite{CDG}.
As we will show here, this limiting, high-dimensional single scale problem can be
solved numerically by sparse tensor product FEM in complexity
comparable to the optimal one for single scale problems in the physical domain
$D$.

We also show how the same idea can be applied to the fast deterministic calculation
of two and $M$ point spatial correlation functions of random solutions to elliptic
problems.

\section{Homogenization problem} \label{section 2}\setzero
\vskip-5mm \hspace{5mm}

In a bounded domain $D \subset \R^d$ with Lipschitz boundary $\Gamma = \partial D$,
we consider the elliptic problem in divergence form
\begin{equation}\label{2.1}
\mbox{$-{\rm div}(A^\e(x) \nabla u^\e) = f(x)$ in $D$, $u^\e = 0$ on $\partial D$}\,.
\end{equation}
Problem (\ref{2.1}) has multiple separated scales in the sense that
\begin{equation}\label{2.2}
A^\e(x) = A\Big(x, \Big\{\frac{x}{\e}\Big\}_Y\Big), \;x \in D
\end{equation}
where $Y = (0,1)^d$ denotes the unit cell and we denote by $[\frac{x}{\e}]_Y$ the unique element
in $\Z^d$ such that $x \in \e ([\frac{x}{\e}]_Y + Y)$ and set
$\{\frac{x}{\e}\}_Y : = \frac{x}{\e} - [\frac{x}{\e}]_Y \in Y$.
In (\ref{2.2}), the function $A(x,y) \in L^\infty(D \times Y)^{d \times d}_{\rm sym}$
is $Y$-periodic with respect to $y$ and satisfies, for every $(x,y) \in D \times Y$,
and some $0< \alpha < 1$:
\begin{equation}\label{2.3}
\forall \xi \in \R^d: \alpha \,|\xi|^2 \le \xi^\top\,A(x,y) \,\xi \le \alpha^{-1} |\xi|^2\,.
\end{equation}
Then problem (\ref{2.1}) admits, for every $f \in L^2(D)$, a unique weak solution:
\begin{equation}\label{2.4}
u^\e \in H^1_0(D): \quad
B^\e(u^\e,v) = (f,v) \quad \forall v \in H^1_0(D)
\end{equation}
where $B^\e(u,v) = \int_D \nabla v \cdot A(x,\frac{x}{\e}) \nabla u \,dx$.
As it is well known, as $\e \rightarrow 0$, $u^\e \rightarrow u^0$ in $L^2(D)$ strongly
and in $H^1(D)$ weakly, and $u^0$ is the solution of the homogenized problem
\begin{equation}\label{2.5}
\mbox{$-{\rm div}(A^0(x) \nabla u^0) = f$ in $D$, $u^0 = 0$ on $\partial D$}\,,
\end{equation}
and formulas for $A^0(x)$ are available \cite{BLP78,JKO94}.
The lack of strong $H^1(D)$-convergence indicates that in the limit $\e \rightarrow 0$, information
on the fine scale of $u^\e$ is lost. It can be recovered by calculating so-called correctors
by either differentiating $u^\e$ or by solving a second problem of the type
(\ref{2.1}). Both approaches are not attractive as basis for numerical solution
methods: differentiating a numerical approximation of $u^0$ will reduce convergence rates,
and solving (\ref{2.1}) for correctors amounts to solving a problem akin to the original
one.

A variant of the two-scale convergence, originally due to
\cite{Nguetseng,Allaire}, and recently developed in \cite{CDG},
allows to obtain a single scale, ``unfolded'', limit problem which gives $u^0$ and
essential information on oscillations of $u^\e$.
To describe it, define for every $\varphi \in L^2(D)$ the ``unfolding'' operator
\begin{equation}\label{2.6}
{\cal T}_\e (\varphi)(x,y) : =
\left\{
\begin{array}{ll}
\varphi\Big(\e\Big[\dis\frac{x}{\e}\Big]_Y + \e y\Big) &
\mbox{if $\; \e\Big(\Big[\dis\frac{x}{\e}\Big]_Y + Y\Big) \subset D$}
\\[2ex]
0 & \mbox{else}
\end{array}\right.
\end{equation}
for $x \in D, y \in Y$. Then there holds \cite{CDG}.

\bigskip\noindent
{\bf Proposition 2.1.} {\it
Assume that
\begin{equation}\label{2.7a}
\tilde{A}^\e(x,y)
:=
{\cal T}_{\e}(A^{\e})(x,y) \rightarrow \tilde{A}(x,y) \quad\mbox{a.e.} (x,y)\in D\times Y.
\end{equation}
Then there exists in $u^0 \in H^1_0(D)$ such that, as $\e \rightarrow 0$,
the solutions $u^\e$ of (\ref{2.4}) satisfy
\begin{equation}\label{2.7}
u^\e \underset{H^1(D)}{\mbox{\LARGE $\rightharpoonup$}}\; u^0 \in H^1_0(D)
\end{equation}
and there exists $\phi(x,y) \in L^2(D,H^1_{\rm per}(Y)/\R)$ such that, as $\e \rightarrow 0$,
\begin{equation}\label{2.8}
{\cal T}_\e(\nabla_x u^\e) \underset{L^2(D\times Y)^d}{\mbox{\LARGE $\rightharpoonup$}}\;
\nabla_x u^0 + \nabla_y \phi \;\,{\rm in} \;\,L^2(D \times Y)^d \,.
\end{equation}
The functions $u^0 \in H^1_0(D)$, $\phi \in L^2(D,H^1_{\rm per}(Y)/\R)$ solve the
``unfolded'' limiting problem: find
$u^0 \in H^1_0(D)$, $\phi \in L^2(D,H^1_{\rm per}(Y)/\R)$ such that
\begin{equation}\label{2.9}
\begin{array}{l}
B(u^0,\phi;v,\psi) : =
\\[1ex]
\dis\int\limits_D \;\dis\int\limits_{Y}
\;\{\nabla_x v + \nabla_y \psi\}^\top \,\tilde{A}(x,y) \{\nabla_x u^0 + \nabla_y \phi\} dy\,dx =
\\[1ex]
\dis\int\limits_D \,fv \,dx \quad \forall v \in H^1_0(D), \,\psi \in L^2(D,H^1_{\rm per} (Y)/\R) \,.
\end{array}
\end{equation}
}

\noindent
{\bf Remark 2.2.}
i) Problem (\ref{2.9}) is independent of $\e$ and therefore a single scale problem.

\medskip\noindent
ii) As $\e \rightarrow 0$, $u^\e(x) \rightarrow u^0(x)$ in $L^2(D)$ and
$\nabla_x u^\e \rightarrow \nabla_x u^0 + \nabla_y\phi(\cdot,\frac{\cdot}{\e})$
strongly in $L^2(D,\R^d)$.
Therefore, $u^0$ coincides with the solution of the homogenized problem (\ref{2.5}) and
$\phi$ retains information on the oscillations of $u^\e$ as $\e \rightarrow 0$.
Information on gradients of $u^\e$ is introduced
at the price of higher dimension of the limiting problem:
while (\ref{2.5}) is posed on $D \subset \R^d$,
(\ref{2.9}) must be solved on $D \times Y \subset \R^{2d}$.

\medskip\noindent
iii) As (\ref{2.9}) is derived by weak convergence methods, no information about the boundary behaviour of $u^\e$
is obtained by solving (\ref{2.9}).

\medskip\noindent
iv) In (\ref{2.4}), $f = f(x)$ was assumed independent of $\e$. Loadings of the type $f(x,x/\e)$ may equally be accommodated.

\medskip \noindent
v) Property (\ref{2.7a}) holds under certain conditions on $A^\e$, for example if
$A(x,y) = A(y)$ or if $A\in L^1(Y,C^0(D))$.

\bigskip
The limit problem (\ref{2.9}) is well-posed:

\medskip\noindent
{\bf Proposition 2.3.} {\it Assume (\ref{2.3}). Then there is $C >
0$ such that
\begin{equation}\label{2.10}
\begin{array}{l}
\forall u \in H^1_0(D), \;\forall \phi \in L^2(D,H^1_{\rm per}(Y)/\R):
\\[2ex]
B(u,\phi; u,\phi) \ge C\big(\|u\|^2_{H^1(D)} + \|\phi\|^2_{L^2(D,H^1_{\rm per}(Y))}) \,.
\end{array}
\end{equation}
}

The proof is immediate if we observe (\ref{2.3}) and
\begin{equation*}
|\nabla_x u + \nabla_y \phi|^2 = |\nabla_x u|^2 + 2 \nabla_x u \cdot \nabla_y \phi + |\nabla_y \phi|^2
\end{equation*}
and that, due to $\phi \in L^2(D,H^1_{\rm per} (Y)/\R)$, it holds
\begin{equation*}
\dis\int\limits_D \;\dis\int\limits_Y \;\nabla_x u \cdot \nabla_y \phi \,dy\,dx = \dis\int\limits_D \,\nabla_x u \cdot \dis\int\limits_Y \,\nabla_y \,\phi\,dy \,dx = 0 \,.
\end{equation*}
The limit problem (\ref{2.9}) admits the following regularity: if $\tilde{A}(x,y)$, $\partial D$
are smooth, then $f \in H^{-1 + k}(D)$ implies
\begin{equation}\label{2.11}
u^0 \in H^{1+k}(D), \;\phi \in H^k(D,C^\infty_{\rm per}(\overline{Y}))\,.
\end{equation}

\section{Sparse finite element discretization} \label{section 3}\setzero
\vskip-5mm \hspace{5mm}

We discretize the unfolded limiting problem (\ref{2.9}) by a sparse Finite Element Method (FEM) (e.g. \cite{GOS} and the references there). To this end, assume $D$ is a bounded Lipschitz polyhedron with straight faces and
let $\{{\cal T}^D_\ell\}$, $\{{\cal T}_\ell^Y\}$ be sequences of nested, quasiuniform meshes in $D$
resp. $Y$ consisting
of shape-regular simplices $T$ of meshwidth $h_\ell = O(2^{-\ell})$, and such that
the periodic extension of ${\cal T}^Y_\ell$ beyond $Y$ is regular.

 Let further $p \ge 1$ be a polynomial degree. Then we denote
\begin{equation*}
\begin{split}
V^{\ell,1}_D & = S^{p,1}_0(D, \cT^D_\ell) = \{u \in H^1_0(D): \forall T \in \cT_\ell^D: u|_T \in {\cal P}_p(T)\}\,,
\\[2ex]
V^{\ell}_Y & = S^{p,1}_{\rm per} (Y,\cT^Y_\ell) = \{u \in H^1_{\rm per}(Y): \forall T \in \cT_\ell^Y: u|_T \in {\cal P}_p(T)\} / \R\,,
\\[2ex]
V^{\ell,0}_D & = S^{p-1,0}(D,\cT^D_\ell) = \{u \in L^2(D): \forall T \in \cT_\ell^D: u|_T \in {\cal P}_{p-1}(T)\}\,,
\end{split}
\end{equation*}
where ${\cal P}_p(T)$ denotes the polynomials of total degree at most $p$ on $T$.

\medskip
Since the triangulations are nested, the Finite Element spaces are hierarchical:
\begin{equation}\label{3.1}
\begin{array}{lllllll}
V^{0,j}_D &\subset &V^{1,j}_D &\subset \dots \subset &V^{\ell,j}_D &\subset \dots \subset &H^j(D) \,,
\\[1ex]
V^0_Y &\subset &V^1_Y &\subset \dots \subset &V^\ell_Y &\subset \dots \subset& H^1_{\rm per}(Y)\,.
\end{array}
\end{equation}
Define orthogonal projections $P^{\ell,j}_D: H^j(D) \rightarrow V^{\ell,j}_D$, $j = 0,1$ and $P^\ell_Y: H^1_{\rm per}(Y) \rightarrow V^\ell_Y$ and set $P_D^{-1,j} : = 0$, $P_Y^{-1} : = 0$.
Then we have the increment or detail-spaces
\begin{equation}\label{3.2}
\begin{array}{llll}
W_D^{\ell,j} & : = &(P^{\ell,j}_D - P_D^{\ell-1,j}) \,V^{\ell,j}_D, &\ell = 0,1,2,\dots
\\[1ex]
W^\ell_Y & : = &(P^\ell_Y - P_Y^{\ell - 1}) \,V^\ell_Y, &\ell = 0,1,2,\dots \end{array}
\end{equation}
and, for every $L > 0$, the multilevel decompositions
\begin{equation}\label{3.3}
V^{L,j}_D = \bigoplus\limits_{0 \le \ell \le L} W_D^{\ell,j},
\qquad
V^L_Y = \bigoplus\limits_{0 \le \ell \le L} \,W_Y^\ell, \quad j = 0,1
\end{equation}
which are orthogonal in $H^j(D)$ resp. $H^1_{\rm per}(Y)$.

\medskip
For $L > 0$, define the full tensor product space
\begin{equation}\label{3.4}
S^L = V^{L,0}_D \otimes V^L_Y = \bigoplus\limits_{0 \le \ell, \ell^\prime \le L} \,W^{\ell,0}_D \otimes W^{\ell^\prime}_Y \subset L^2(D,H^1_{\rm per}(Y)/\R) \,.
\end{equation}
Assume the regularity (\ref{2.11}) with some $k > 0$.

\medskip
Then (\ref{2.10}) shows that the finite element approximation
\begin{equation}\label{3.5}
(u^L,\phi^L) \in V^{L,1}_D \times S^L : B(u^L,\phi^L; v^L,\psi^L) = (f,v^L) \quad \forall (v^L,\psi^L) \in V^{L,1}_D \times S^L
\end{equation}
exists, is unique and a tensor product argument show that it satisfies the asymptotic error bounds
\begin{equation}\label{3.6}
\begin{array}{l}
\|u^0 - u^L\|_{H^1(D)} + \|\phi - \phi^L \|_{L^2(D,H^1(Y))}
\\[1ex]
\qquad \le C \,h_L^{\min(p,k)} \,\big\{\|u^0\|_{H^{k+1}(D)} + \|\phi\|_{L^2(D; H^{k+1}(Y))} + \|\phi\|_{H^k(D;H^1(Y))}\big\} \,.  \end{array}
\end{equation}
We note that dim$(V^{L,1}_D) = O(2^{L d}) = {\rm dim}\,V^L_Y$, whereas dim$(S^L) = O(2^{2 L d}) = O({\rm dim}(V^{L,1}_D)^2)$ as $L \rightarrow \infty$, since $S^L$ is a FE-space in $D \times Y \subset \R^{2d}$. Also, note that in (\ref{3.6}) the full regularity (\ref{2.11}) of $\phi$ was not used.

The large number of degrees of freedom in $S^L$ due to high dimension renders the unfolded problem (\ref{2.9})
impractical for efficient numerical solution.
A remedy is to use the sparse tensor product spaces (\cite{GOS} and the references there)
\begin{equation}\label{3.7}
\wh{S}^L : = \bigoplus\limits_{0 \le \ell + \ell^\prime \le L}\,W^{\ell,0}_D \otimes W^{\ell^\prime}_Y
\end{equation}
which have substantially smaller dimension than the full tensor product spaces
(\ref{3.4}).

The joint regularity (\ref{2.11}) in $x$ and $y$ of $\phi(x,y)$ allows to retain the
error bounds (\ref{3.6}) with $\wh{S}^L$ in place of $S^L$:

\medskip\noindent
{\bf Proposition 3.1.}
{\it
Let $(\wh{u}^L, \wh{\phi}^L) \in V^{L,1}_D \times \wh{S}^L$ denote the sparse FE-solutions
which satisfy (\ref{3.5}) with $\hat{S}^L$ in place of $S^L$.

If $u^0 \in H^{1 + k}(D)$, $\phi \in H^k(D,H^{k+1}_{\rm per} (Y))$, for some $k \ge 1$, then
\begin{equation}\label{3.8}
\begin{array}{l}
\|u^0 - \wh{u}^L\|_{H^1(D)} + \|\phi - \wh{\phi}^L \|_{L^2(D,H^1(Y))}
\\[1ex]
\qquad \le C  (L+1)^{\frac{1}{2}} \,h_L^{\min(p,k)} \,\big\{\|u^0\|_{H^{k+1}(D)} + \|\phi\|_{H^k(D; H^{k+1}(Y))}
\big\}
\end{array}
\end{equation}
and the total number of degrees of freedom is bounded by
\begin{equation}\label{3.9}
{\rm dim}(V^{L,1}_D) + {\rm dim} (\wh{S}^L) \le C L 2^{dL} \,.
\end{equation}
}

\noindent
{\bf Remark 3.2.}
i) Since dim$(V^{L,1}_D) \le C 2^{dL}$, computation of $(\wh{u}, \wh{\phi}^L)$ in $D \times Y$ requires,
up to $L = O(|\log h_L|)$, the same number of degrees of freedom as the FE approximation of
$u^0$ from (\ref{2.5}) in $D$.
Moreover, the FE approximation of (\ref{2.9}) does not require the determination of $A^0(x)$.
In addition, the convergence rate (\ref{3.8}) is, up to $(L+1)^{\frac{1}{2}}$,
equal to the rate (\ref{3.6}) for the full tensor product-spaces.

\medskip\n
ii) Sparse tensor products of
finite elements in one dimension were proposed by Zenger and his students
in the 1990ies for the efficient solution of partial differential equations in
three dimensions (see \cite{GOS} and the references there).
While allowing similar convergence rates as the full tensor product spaces,
extra regularity of the solution (generally not available in non-smooth domains) is required.
In (\ref{3.7}), sparse tensor products
of standard finite element spaces in $D$ resp. $Y\subset \R^d$ are taken and
realistic regularity of the solution available from the structure of the limit problem
(\ref{2.9}) was used.

\medskip\n
iii) The approach generalizes to problems with $M > 2$ scales.
The unfolded problem is then posed on a product domain
\begin{equation*}
D \times Y_1 \times \dots \times Y_{M-1} \subset \R^{Md}
\end{equation*}

\medskip\n
and an approximation of order $(L+1)^{\frac{M-1}{2}} h_L^{\min(p,k)}$ can be obtained with
$O(L^{M-1} \,2^{dL})$ degrees of freedom.

\medskip\n
iv) Explicit construction of the sparse tensor product-space $\wh{S}^L$ requires bases for the
detail-spaces $W^{\ell,0}_D$, $W^\ell_Y$. These are available, for example, via suitable semiorthogonal
wavelet bases (e.g. \cite{Dah97}). These bases allow also for optimal preconditioning
of the linear system corresponding to (\ref{3.5}).

\medskip\n
v) We discussed here only the diffusion problem (\ref{2.1}). The results can be generalized to Elasticity,
and the Stokes Equations.

\medskip\n
vi)
If only partial periodicity or patchwise periodic patterns are present, the unfolding approach works equally.
See \cite{CDG} for details.

\section{Stochastic data} \label{section 4}\setzero
\vskip-5mm \hspace{5mm}

Let $(\Omega,\Sigma,P)$ be a $\sigma$-finite probability space and $A(x) \in L^\infty(D,\R^{d \times d}_{\rm sym})$ satisfy for every $\xi \in \R^d$
\begin{equation}\label{4.1}
\exists \alpha, \beta > 0: \alpha |\xi|^2 \le \xi^T \,A(x) \,\xi \le \beta |\xi|^2 \,.
\end{equation}

\n
For a random source term $f \in L^2(\Omega,d P; L^2(D))$, consider the Dirichlet
problem \begin{equation}\label{4.2}
\mbox{$L(x,\partial_x) u = - \nabla \cdot A(x) \,\nabla u(x,\omega)
=
f(x,\omega)$ in $D$, $u = 0$ on $\partial D$}
\end{equation}
for $P$-a.e. $\omega \in \Omega$. The random solution $u(x,\omega)$ of (\ref{4.2})
is searched in the Bochner-space
\begin{equation}\label{4.3}
{\cal H}^1_0(D) : = L^2(\Omega, d P; H^1_0(D)) \cong L^2(\Omega, dP) \otimes H^1_0(D) \,.
\end{equation}

\smallskip\n
We note that $L^2(\Omega,dP)$, equipped with inner product
\begin{equation}\label{4.4}
\langle u,v\rangle = \dis\int\limits_\Omega \,u(\omega) \,v(\omega) \,dP(\omega) \,,
\end{equation}
is a Hilbert space. The variational form of (\ref{4.2}) reads: find $u \in {\cal H}^1_0(D)$ such that
\begin{equation}\label{4.5}
{\cal A}(u,v) = \langle f,v \rangle_{{\cal L}^2(D)} \quad \forall v \in {\cal H}^1_0(D)
\end{equation}
where
\begin{equation*}
{\cal A}(u,v) = \langle (A \otimes id) (\nabla \otimes id) u, \;(\nabla \otimes id) v \rangle_{{\cal L}^2(D)^d} \,
\end{equation*}

\medskip\n
and $\langle f,v \rangle_{{\cal L}^2(D)} = \int_D \,\langle f(x,\cdot), v(x, \cdot)\rangle dx$.
The form ${\cal A}(\cdot,\cdot)$ in (\ref{4.5}) is coercive on ${\cal H}^1_0(D) \times {\cal H}^1_0(D)$,
implying the existence of a unique random solution $u(x,\omega)$ of (\ref{4.5}).
Numerical solution of (\ref{4.5}) that involves a FEM in $D$ and Monte-Carlo in $\Omega$ is prohibitively
expensive. Alternatively, we might try to compute directly the statistics of $u(x,\omega)$.
For example, the mean field $E_u(x) = \int_\Omega u(x,\omega) dP(\omega)$ solves
\begin{equation}\label{4.6}
\mbox{$L(x,\partial_x) E_u = E_f$ in $D$, $E_u = 0$ on $\partial D$}\,.
\end{equation}
For $x,x^\prime \in D$, define the two point correlation function
\begin{equation}\label{4.7}
C_u(x,x^\prime) = \langle u(x,\cdot), u(x^\prime,\cdot)\rangle \,.
\end{equation}
Then the variance of $u(x,\omega)$ is given by
\begin{equation}\label{4.8}
({\rm Var} \,u(x,\cdot))^2 : = (E_u(x))^2 - (C_u(x,x))^2 \,.
\end{equation}
The two-point correlation $C_u(x,x^\prime)$ is the solution of a deterministic problem in $D \times D$.
Formally
\begin{equation}\label{4.9}
L(x,\partial_x) \,L(x^\prime,\partial_{x^\prime}) C_u = C_f \;\,{\rm in} \;\,D \times D \,,
\end{equation}
and in variational form:
\begin{equation}\label{4.10}
C_u \in H^{1,1}_{(0)}(D \times D): Q (C_u,C_v) = (C_f,C_v) \quad \forall C_v \in H^{1,1}_{(0)}(D\times D)
\end{equation}
where $H^{1,1}_{(0)}(D\times D) := H^1_0(D) \otimes H^1_0(D)$ and
\begin{equation*}
Q(C_u,C_v) = \dis\int\limits_{D \times D} \,\nabla_{xy} \,C_v \;A(x) \otimes A(y) \,\nabla_{xy} \,C_u \,dx\,dy, \;\nabla_{xy} := \nabla_x \otimes \nabla_y \,.
\end{equation*}

\medskip\n
{\bf Proposition 4.1.} {\it The two point correlation $C_u$ of the random solution $u(x,\omega)$ is the unique solution of the deterministic problem (\ref{4.10}) in $D \times D$.}

\medskip
Hence to get two point correlation functions, Monte-Carlo can be traded for a deterministic problem
in high dimensions. The key to its efficient numerical solution lies in the regularity of the solution:
if the mean field problem (\ref{4.6}) satisfies a shift-theorem at order $s > 0$,
i.e. $E_f \in H^{-1+s}(D) \Longrightarrow E_u \in H^{1+s}(D)$, then
\begin{equation}\label{4.11}
C_f \in H^{-1+s,-1+s}(D \times D) \Longrightarrow C_u \in H^{1+s,1+s}_{(0)}(D \times D) \,.
\end{equation}
A similar regularity result in weighted spaces of mixed highest derivatives holds if $D \subset \R^2$
has corners \cite{ST01}.
This regularity in spaces of mixed highest derivatives allows to approximate $C_u$
from the sparse tensor product FE spaces
\begin{equation}\label{4.12}
\wh{V}^L : = \bigoplus\limits_{0 \le \ell + \ell^\prime \le L} \,W^{\ell,1}_0 \otimes W^{\ell^\prime,1}_D
\end{equation}
of dimension dim$(\wh{V}^L) \le C L 2^{dL}$ at a near optimal convergence rate: if $\wh{C}^L_u \in \wh{V}^L$ denotes the FE approximation of $C_u$, it holds \cite{ST01}
\begin{equation}\label{4.13}
\|C_u - \wh{C}^L_u\|_{H^{1,1}(D \times D)} \le C \,|\log h_L|^{\frac{1}{2}}\,h_L^{\min(s,p)} \|C_u\|_{H^{s+1,s+1}(D \times D)}\,.
\end{equation}
Moreover, using a semiorthogonal wavelet basis of the detail-spaces $W^{\ell,1}_D$ in (\ref{3.2}), we can design an algorithm which computes $\wh{C}^L_u$ to the order of the discretization error (\ref{4.13}) in $O(N_L L^{4d+2})$ operations where $N_L = O(2^{dL})$ denotes the number of degrees of freedom in $D$ (see \cite{ST01} for details).

\medskip
As for homogenization problems with multiple scales, $M$ point correlation functions of $u(x,\omega)$ can be approximated at the rate $|\log h_L|^{(M-1)/2}$ $h_L^{\min(s,p)}$ with $O(N_L \,L^{M-1})$ degrees of freedom . Let us also remark that high regularity in (\ref{4.13}) corresponds to strong spatial correlation of $u(x,\omega)$.

In the spatially uncorrelated limit, formally
$C_f(x,y) = \delta(x-y)$ and we have for smooth $\partial D, A(x),d \le 3$:
\begin{equation*}
C_u \in H^{1 + s,1+s}(D \times D), \;0\leq s < 1 - \dis\frac{d}{4}\,,
\end{equation*}
so that only the low convergence rates $|\log h_L|^{\frac{1}{2}}\,h_L^{1-\frac{d}{4} - \e}$
for $\wh{C}^L_u$ follow.
This is due to the singular support of $C_u$ being the diagonal $\{(x,y) : x = y\}$.
The efficient approximation of such $C_u$ is the topic of ongoing research.

\label{lastpage}

\end{document}